\documentclass[11pt]{amsart}
\usepackage[T1]{fontenc}
\usepackage[utf8]{inputenc}
\usepackage{lmodern}
\usepackage[margin=1.0in]{geometry}
\usepackage{microtype}
\usepackage{amsmath,amssymb,amsthm,mathtools}
\usepackage{enumitem}
\usepackage{xcolor}
\usepackage[backref=page]{hyperref}
\usepackage{mathrsfs}
\usepackage{thmtools}
\usepackage{bookmark}

\hypersetup{colorlinks=true,linkcolor=blue,citecolor=purple,urlcolor=blue}

\renewcommand*{\backref}[1]{}
\renewcommand*{\backrefalt}[4]{%
  {\color{blue}%
    \ifcase #1 %
    \or (#2)%
    \else (#2)%
    \fi}}

\declaretheoremstyle[
  headfont=\bfseries,
  bodyfont=\itshape,
  spaceabove=0.8em,
  spacebelow=0.8em,
  headpunct={.},
  postheadspace=0.5em
]{plainstyle}
\declaretheoremstyle[
  headfont=\bfseries,
  bodyfont=\normalfont,
  spaceabove=0.8em,
  spacebelow=0.8em,
  headpunct={.},
  postheadspace=0.5em
]{defstyle}
\declaretheorem[style=plainstyle,numberwithin=section]{theorem}
\declaretheorem[style=plainstyle,sibling=theorem]{proposition}

\declaretheorem[style=plainstyle,sibling=theorem]{corollary}
\declaretheorem[style=defstyle,sibling=theorem]{definition}
\declaretheorem[style=defstyle,sibling=theorem]{remark}
\declaretheorem[style=defstyle,sibling=theorem]{example}



\newcommand{\cml}{\mathrm{cml}}
\newcommand{\Ker}{\operatorname{ker}}

\newcommand{\gen}[1]{\langle #1 \rangle}

\newcommand{\FF}{\mathbb{F}}

\newcommand{\betafun}{\beta}

\title{Approximate Subloops in Moufang Loops}
\author{Arindam Biswas}
\email{arin.math@gmail.com}
\date{\today}

\begin{document}

\begin{abstract}
We introduce a notion of finite approximate subloops in Moufang loops, with emphasis on the commutative case. For arbitrary Moufang loops we establish intrinsic product-set identities and covering consequences without passing through associative quotients and obtain a finite-kernel reduction principle: approximate-subloop structure descends through homomorphisms onto groups with finite kernel, and inverse results in the quotient lift back to the loop. In particular, this yields a complete reduction in the two-generated case. For commutative Moufang loops, using their local finite-by-abelian structure, we deduce a Freiman-type theorem showing that a finite approximate subloop is contained in the pullback of a coset progression from a suitable local abelian quotient, with quantitative bounds depending only on the corresponding finite kernel. We then obtain a uniform version for approximate subloops generating an $m$-generated subloop. When the local abelian quotient has bounded torsion, we get a polynomial covering theorem by cosets of a finite subloop, deduced from the bounded-torsion polynomial Freiman--Ruzsa theorem in the abelian quotient; in particular, this applies to commutative Moufang loops of exponent $3$.
\end{abstract}

\maketitle

\section{Introduction}

Approximate groups and small-doubling sets are central objects in additive and
multiplicative combinatorics. In abelian groups, Freiman-type inverse theorems
describe finite sets of small doubling in terms of coset progressions
\cite{Freiman64,GR07}. In nonabelian groups, the modern theory of approximate
groups has developed through work of Tao, Hrushovski, Breuillard--Green--Tao,
Pyber--Szab\'o, Helfgott, and others
\cite{Tao08,Hrushovski12,BGT11,BGT12,PyberSzabo16}. The purpose of the present
paper is to formulate and prove loop-theoretic analogues of some of these
structural results for Moufang loops, with particular emphasis on the
commutative Moufang setting.

A basic difficulty in the nonassociative setting is that higher product sets are
not canonically defined: while a twofold product $AB$ is unambiguous, triple and
higher products depend on bracketing. Thus one must decide whether to work with
a fixed bracketing convention or with an intrinsic notion that records all
possible bracketings. In this paper we adopt the latter point of view. This is
one of the reasons that we use a small-tripling definition of approximate
subloop rather than a covering definition built from a preferred product.

A second structural issue is that, unlike in the group case, one cannot in
general pass immediately to an associative ambient object. The relevant
associative information in a Moufang loop often appears only after quotienting by
a suitable finite kernel. In the commutative Moufang setting this becomes
especially effective: if $A$ is finite and $H=\gen{A}$, then the associator
subloop $H'=(H,H,H)$ is finite and the quotient $H/H'$ is an abelian group
\cite{Bruck58,Pflugfelder90}. This local finite-by-abelian structure allows one
to combine intrinsic loop-theoretic arguments in $H$ with additive-combinatorial
inverse theorems in the abelian quotient.

The results of the paper fall into three parts.. First, for arbitrary Moufang loops, we prove intrinsic product-set identities and covering consequences that do not pass through associative quotients. We
also establish a finite-kernel reduction principle: if a Moufang loop admits a
homomorphism onto a group with finite kernel, then approximate-subloop
structure descends to the quotient, and inverse results in the quotient lift
back to the loop. In particular, whenever a finite set generates a
two-generated Moufang subloop, the problem is entirely associative, since every
two-generated Moufang loop is a group \cite{Bruck58,Pflugfelder90}. Second, for commutative Moufang loops, we use the local quotient
\[
H/H', \qquad H=\gen{A}, \quad H'=(H,H,H),
\]
to obtain a local Freiman-type structure theorem: a finite approximate subloop
is contained in the pullback of a coset progression from the abelian quotient,
with quantitative bounds depending only on the finite kernel $H'$. We then give
a uniform version when the generated subloop is $m$-generated. Third, when the local abelian quotient has bounded torsion, we obtain a polynomial covering theorem by cosets of a finite subloop, deduced from the bounded-torsion polynomial Freiman--Ruzsa theorem of
Gowers--Green--Manners--Tao \cite{GGMT26}. In particular, this applies to
commutative Moufang loops of exponent $3$.

In Section~\ref{sec:preliminaries} we collect the background and
notation used throughout, including our conventions for associators and
bracketed product sets. Section~\ref{sec:approximate-subloops} introduces approximate subloops and develops
the intrinsic product-set identities and covering lemmas used later. Section~\ref{sec:finite-kernel}
contains the finite-kernel reduction principle in arbitrary Moufang loops.
Sections~\ref{sec:pullbacks} and~\ref{sec:local-structure} specialize to commutative Moufang loops, proving the local and
uniform structure theorems. Section~\ref{sec:bounded-torsion} treats the bounded-torsion case and the
exponent-$3$ application.

\section{Preliminaries and notation}
\label{sec:preliminaries}

\subsection{Loops and Moufang loops}

\begin{definition}[Quasigroups and loops]
A \emph{quasigroup} is a nonempty set $Q$ with a binary operation
$(x,y)\mapsto xy$ such that for every $a,b\in Q$ the equations
\[
ax=b, \qquad ya=b
\]
have unique solutions $x,y\in Q$. A \emph{loop} is a quasigroup with a
two-sided identity element $1$, that is,
\[
1x=x1=x \qquad (x\in Q).
\]
\end{definition}

For a loop $L$ and $a\in L$, we denote the left and right translations by
\[
L_a(x)=ax, \qquad R_a(x)=xa.
\]

\begin{definition}[Moufang loop \cite{Mou35}]
A loop $L$ is a \emph{Moufang loop} if it satisfies
\[
(xy)(zx)=(x(yz))x \qquad (x,y,z\in L).
\]
A Moufang loop is \emph{commutative} if $xy=yx$ for all $x,y\in L$.
\end{definition}

We use the standard consequences of the Moufang
identities: Moufang loops are alternative, they satisfy the inverse
property, inversion is an anti-automorphism, and every two elements generate a
group. In particular, powers are unambiguously defined, and any expression
involving only two fixed elements and their powers may be rebracketed freely
\cite{Bruck58,Pflugfelder90}.

\begin{example}
Every group is a Moufang loop. A classical nonassociative commutative Moufang
loop is the Zassenhaus--Bol example on $\FF_3^4$ with multiplication
\[
x\circ y=x+y+\bigl(0,0,0,(x_3-y_3)(x_1y_2-x_2y_1)\bigr).
\]
This gives a concrete commutative Moufang loop that is not associative.
Historically, nonassociative commutative Moufang loops go back to Bol's 1937
paper, where the construction is attributed to Zassenhaus \cite{Bol37}; for the
displayed coordinate model above, see also \cite{Smith16}.
\end{example}

\subsection{Associators and quotient loops}

\begin{definition}[Associators]
Let $L$ be a loop. The \emph{associator} $(x,y,z)$ of $x,y,z\in L$ is defined by
\[
(xy)z=x(yz)\,(x,y,z).
\]
The \emph{associator subloop} of $L$ is the subloop generated by all
associators:
\[
(L,L,L):=\gen{(x,y,z):x,y,z\in L}.
\]
\end{definition}

If $N\le L$ is a normal subloop, we write $L/N$ for the corresponding quotient
loop and $\pi:L\to L/N$ for the quotient map. When $N$ is finite, every fiber
of $\pi$ is a coset of $N$ and has cardinality $|N|$.

The following classical facts about commutative Moufang loops will be used
repeatedly.

\begin{theorem}[Classical structure of commutative Moufang loops]
\label{thm:cml-classical}
Let $M$ be a commutative Moufang loop.
\begin{enumerate}[label=(\roman*)]
  \item The associator subloop $M':=(M,M,M)$ is normal.
  \item The quotient $M/M'$ is an abelian group.
  \item The cube map $x\mapsto x^3$ is an endomorphism of $M$ with image
  contained in $Z(M)$; in particular $x^3\in Z(M)$ for all $x\in M$.
  \item If $M$ is finitely generated, then $M$ is centrally nilpotent.
  \item If $M$ is finitely generated, then $M'$ is finite.
\end{enumerate}
\end{theorem}

\begin{proof}
These are classical structural properties of commutative Moufang loops; see
\cite{Bruck58,Pflugfelder90}. For the associator calculus used in finitely
generated commutative Moufang loops, see also \cite{Smith78}.
\end{proof}

\begin{remark}[The Triple Argument Hypothesis]
In a loop, the \emph{Triple Argument Hypothesis} (TAH) is the assertion that
every associator with one argument appearing three times is trivial. In the
commutative Moufang setting, Smith recorded the problem of whether the TAH
follows from the commutative Moufang identities \cite{Smith78}.
This is now known to be false: Sandu refuted the conjecture for commutative
Moufang loops, and Grishkov--Shestakov later showed that it already fails in
the free commutative Moufang loop of exponent $3$ on $7$ generators
\cite{Sandu88,GrishkovShestakov11}. We do not assume the TAH anywhere.
\end{remark}

\subsection{Subset products and bracketing conventions}

For subsets $X,Y$ of a loop $L$, we write
\[
XY:=\{xy:x\in X,\ y\in Y\}.
\]
For two subsets there is no ambiguity, but for triple and higher products one
must specify bracketings.

\begin{definition}
Let $A$ be a finite subset of a loop $L$.
\begin{enumerate}[label=(\roman*)]
  \item We write $A^2:=AA$.
  \item We write
  \[
  A^{\langle 3\rangle}:=(AA)A\cup A(AA).
  \]
  \item More generally, for $n\ge 2$, $A^{\langle n\rangle}$ denotes the union of
  all subset products obtained from $n$ copies of $A$ using all possible
  bracketings. The number of such subsets will be $C_{n-1}$ the nth Catalan number.
\end{enumerate}
\end{definition}

For example, when $n=4$ there are five bracketings,
\[
((AA)A)A,\qquad (A(AA))A,\qquad (AA)(AA),\qquad A((AA)A),\qquad A(A(AA)),
\]
and $A^{\langle 4\rangle}$ is the union of the corresponding five subset
products. In an associative setting one recovers the usual notation
$A^{\langle n\rangle}=A^n$.

If $L$ is a commutative Moufang loop, then passage to the quotient
$L/(L,L,L)$ removes all bracketing ambiguity, since
Theorem~\ref{thm:cml-classical}\textnormal{(i)}--\textnormal{(ii)} shows that this
quotient is associative and commutative. This is the basic link between
nonassociative product sets in commutative Moufang loops and classical
additive-combinatorial arguments in abelian groups.

\section{Approximate subloops}
\label{sec:approximate-subloops}
We adopt a small-tripling definition. It is stable under quotienting and does not require a preferred higher-product bracketing.

\begin{definition}
Let $L$ be a loop and let $K\ge 1$. A finite subset $A\subseteq L$ is called a \emph{$K$-approximate subloop} if:
\begin{enumerate}[label=(\roman*)]
  \item $1\in A$;
  \item $A=A^{-1}$;
  \item $|A^{\langle 3\rangle}|\le K|A|$.
\end{enumerate}
\end{definition}

\begin{remark}
In the associative case this is polynomially equivalent to the usual definitions of approximate subgroup~\cite{BGT12}. When we speak below of an \emph{$L$-approximate subgroup} in a group, we shall mean a finite symmetric set $B$ with $1\in B$ and $|B^3|\le L|B|$.
\end{remark}

\begin{remark}
The use of $A^{\langle 3\rangle}$ rather than a covering definition is deliberate. In a general loop there is no canonical $n$-fold product set unless a bracketing convention is fixed, whereas $A^{\langle 3\rangle}$ is intrinsic and behaves well under homomorphisms.
\end{remark}

\subsection{Intrinsic product-set identities in Moufang loops}
The Moufang identities yield several useful identities for product sets that do not rely on passage to an associative quotient. We note them here for later use.

\begin{proposition}[Moufang shift identities]
Let $L$ be a Moufang loop, let $z\in L$, and let $X,Y\subseteq L$ be finite nonempty subsets. Let $zXz$ denote the (well-defined) bracketed set $\{zxz : x \in X\}$. Then
\[
z\bigl(X(zY)\bigr) = (zXz)Y 
\qquad \text{and} \qquad 
\bigl((Xz)Y\bigr)z = X(zYz).
\]
Consequently,
\[
|X(zY)| = |(zXz)Y| 
\qquad \text{and} \qquad 
|(Xz)Y| = |X(zYz)|.
\]
\end{proposition}

\begin{proof}
The left Moufang identity states that $z(x(zy)) = ((zx)z)y$ for all $x,y,z\in L$. By diassociativity, $(zx)z = z(xz) = zxz$, allowing us to write the right hand side as $(zxz)y$.
Applying this identity pointwise to the subsets $X$ and $Y$, we obtain
\[
z\bigl(X(zY)\bigr) = (zXz)Y.
\]
Since left translation by $z$ is a bijection of the loop $L$, the set on the left-hand side has exactly the same cardinality as $X(zY)$. Therefore,
\[
|X(zY)| = |(zXz)Y|.
\]
The second identity follows symmetrically from the right Moufang identity $((xz)y)z = x(z(yz)) = x(zyz)$ and the fact that right translation by $z$ is also a bijection.
\end{proof}

\begin{proposition}[Mixed left-right translation identity]\label{prop:moufang-transport}
Let $L$ be a Moufang loop, let $x\in L$, and let $X,Y\subseteq L$. Then
\[
(xX)(Yx)=x(XY)x.
\]
Consequently, if $X$ and $Y$ are finite, then
\[
|(xX)(Yx)|=|XY|.
\]
If, in addition, $L$ is commutative, then
\[
(xX)(xY)=x(XY)x,
\]
and hence
\[
|(xX)(xY)|=|XY|.
\]
\end{proposition}

\begin{proof}
Let $u\in X$ and $v\in Y$. By the Moufang identity,
\[
(xu)(vx)=x(uv)x.
\]
Therefore every element of $(xX)(Yx)$ belongs to $x(XY)x$, and conversely every element of $x(XY)x$ has the form $(xu)(vx)$ for suitable $u\in X$, $v\in Y$. Hence
\[
(xX)(Yx)=x(XY)x.
\]
Since left and right translations in a loop are bijections, the map
\[
t\mapsto xtx
\]
is a bijection of $L$, so
\[
|x(XY)x|=|XY|.
\]
This proves the first assertion.

If $L$ is commutative, then $Yx=xY$, so
\[
(xX)(xY)=x(XY)x.
\]
The cardinality statement follows as above.
\end{proof}

\begin{remark}
Even though left translation need not preserve products in a nonassociative loop, the proposition shows that in a Moufang loop the mixed product of a left $x$-translate and a right $x$-translate has exactly the same size as the original product set.
\end{remark}

\begin{proposition}[Translation invariance of certain product-set sizes]
Let $L$ be a Moufang loop, let $x\in L$, and let $A,B\subseteq L$ be finite nonempty subsets. Then
\[
|(xA)B| = |A(Bx^{-1})|
\qquad\text{and}\qquad
|A(Bx)| = |(x^{-1}A)B|.
\]
\end{proposition}

\begin{proof}
For the first identity, apply Proposition~\ref{prop:moufang-transport} with
$X=A$ and $Y=Bx^{-1}$. Since $(Bx^{-1})x=B$, we obtain
\[
(xA)B=x\bigl(A(Bx^{-1})\bigr)x.
\]
Hence
\[
|(xA)B|=|A(Bx^{-1})|.
\]

For the second identity, apply the same proposition with $X=x^{-1}A$ and $Y=B$.
Since $x(x^{-1}A)=A$, we obtain
\[
A(Bx)=x\bigl((x^{-1}A)B\bigr)x,
\]
and therefore
\[
|A(Bx)|=|(x^{-1}A)B|.
\]
\end{proof}

\begin{proposition}[Symmetry of products and bracketings]\label{prop:bracketsproducts}
Let $L$ be a Moufang loop, and let $A,B\subseteq L$ be finite subsets.
\begin{enumerate}[label=(\roman*)]
  \item $|AB| = |B^{-1}A^{-1}|$.
  \item If $A=A^{-1}$ and $B=B^{-1}$, then $|AB| = |BA|$.
  \item If $A=A^{-1}$, then $|(AA)A| = |A(AA)|$.
\end{enumerate}
\end{proposition}

\begin{proof}
Since Moufang loops are inverse-property loops, the inversion map $x\mapsto x^{-1}$ is an anti-automorphism, so
\[
(xy)^{-1}=y^{-1}x^{-1}\qquad(x,y\in L).
\]
Applying this pointwise to the product set $AB$ gives
\[
(AB)^{-1} = B^{-1}A^{-1}.
\]
Because the inversion map preserves cardinality, we deduce
\[
|AB| = |(AB)^{-1}| = |B^{-1}A^{-1}|.
\]
This proves~(i). For~(ii), if $A=A^{-1}$ and $B=B^{-1}$, applying~(i) yields $|AB| = |B^{-1}A^{-1}| = |BA|$.

For~(iii), applying the anti-automorphism to the threefold product gives
\[
\bigl((AA)A\bigr)^{-1} = A^{-1}(A^{-1}A^{-1}).
\]
Since $A=A^{-1}$, the right-hand side is exactly $A(AA)$. Because inversion is a cardinality-preserving bijection, we conclude that the two bracketings have the same cardinality:
\[
|(AA)A| = |A(AA)|.
\]
\end{proof}

\begin{remark}
Proposition~\ref{prop:bracketsproducts} shows that for symmetric sets the two
threefold bracketings $(AA)A$ and $A(AA)$ always have the same cardinality.
Thus control of one bracketing automatically gives control of the other, without
any commutativity assumption or passage to a quotient.
\end{remark}

The previous result doesn't mean that the cardinality of all higher bracketed products are same. 

\begin{example}[A $4$-fold counterexample in a commutative Moufang loop]
\label{ex:fourfold-counterexample}
Let
\[
M=(\mathbb F_3^4,\circ),
\]
where the multiplication is
\[
x\circ y
=
x+y+\bigl(0,0,0,(x_3-y_3)(x_1y_2-x_2y_1)\bigr),
\]
with all coordinates computed modulo $3$.

This is the standard nonassociative commutative Moufang loop of order $81$.
Its identity element is
\[
0:=(0,0,0,0),
\]
and one checks directly that inversion is given by
\[
x^{-1}=-x.
\]
Indeed,
\[
x\circ(-x)
=
x-x+\bigl(0,0,0,(x_3+x_3)(x_1(-x_2)-x_2(-x_1))\bigr)
=
0.
\]

Now set
\[
u:=(0,0,1,0),\qquad
v:=(0,1,2,0),\qquad
w:=(1,0,0,2),
\]
and note that
\[
u\circ w=(1,0,1,2).
\]
Define
\[
d:=u\circ w=(1,0,1,2),
\]
and
\[
A:=\{0,\pm u,\pm v,\pm w,\pm d\}.
\]
Explicitly,
\[
A=
\{
(0,0,0,0),
(0,0,1,0),(0,0,2,0),
(0,1,2,0),(0,2,1,0),
(1,0,0,2),(2,0,0,1),
(1,0,1,2),(2,0,2,1)
\}.
\]
Since inversion is negation, $A=A^{-1}$.

A direct calculation from the multiplication law gives
\[
AA
=
\{(0,i,j,0):i,j\in\mathbb F_3\}
\cup
\{(1,0,t,2):t\in\mathbb F_3\}
\cup
\{(2,0,t,1):t\in\mathbb F_3\}
\]
\[
\qquad\cup
\{
(1,1,0,1),(1,1,2,0),(1,2,1,0),(1,2,2,2),
(2,1,1,1),(2,1,2,0),(2,2,0,2),(2,2,1,0)
\},
\]
so in particular
\[
|AA|=23.
\]

Multiplying once more by $A$, one finds that $(AA)A=A(AA)$ has complement
\[
M\setminus (AA)A
=
\{(0,0,t,\lambda):t\in\mathbb F_3,\ \lambda\in\{1,2\}\}
\]
\[
\qquad\cup
\{(0,1,2,\lambda),(0,2,1,\lambda):\lambda\in\{1,2\}\}
\]
\[
\qquad\cup
\{(1,0,t,\lambda):t\in\mathbb F_3,\ \lambda\in\{0,1\}\}
\cup
\{(2,0,t,\lambda):t\in\mathbb F_3,\ \lambda\in\{0,2\}\}.
\]
Hence
\[
|(AA)A|=|A(AA)|=81-22=59.
\]

Now compute the two $4$-fold bracketings
\[
((AA)A)A
\qquad\text{and}\qquad
(AA)(AA).
\]
A direct calculation yields
\[
M\setminus ((AA)A)A=\{(0,0,0,1),(0,0,0,2)\},
\]
so
\[
|((AA)A)A|=81-2=79.
\]
On the other hand,
\[
M\setminus (AA)(AA)
=
\{
(0,0,0,1),(0,0,0,2),
(0,1,2,1),(0,1,2,2),
(0,2,1,1),(0,2,1,2)
\},
\]
so
\[
|(AA)(AA)|=81-6=75.
\]

Therefore
\[
|((AA)A)A|=79\neq 75=|(AA)(AA)|.
\]
Thus even for a symmetric set in a commutative Moufang loop, not all
$4$-fold bracketings have the same cardinality.
\end{example}

However certain bracketing preserves the cardinality. This is detailed in the proposition below.

\begin{proposition}[Mirror symmetry for bracketed product sets]
\label{prop:mirror-bracketing}
Let $L$ be a Moufang loop. For each $r\ge 1$, let $T$ be a bracketing of an
$r$-fold product, and let $T^\ast$ denote the mirror bracketing obtained by
interchanging left and right at every binary node. Then for any subsets
$X_1,\dots,X_r\subseteq L$ one has
\[
T(X_1,\dots,X_r)^{-1}
=
T^\ast(X_r^{-1},\dots,X_1^{-1}).
\]
Consequently, if $A=A^{-1}$, then for every bracketing $T$ of an $r$-fold
product of $A$,
\[
|T(A,\dots,A)|=|T^\ast(A,\dots,A)|.
\]
\end{proposition}

\begin{proof}
We first fix the recursive notation.

For $r=1$, there is only one bracketing, and we set
\[
T(X_1)=X_1.
\]
For $r\ge 2$, every bracketing $T$ can be written uniquely in the form
\[
T=(U,V),
\]
where $U$ is a bracketing of the first $k$ variables and $V$ is a bracketing of
the last $r-k$ variables for some $1\le k\le r-1$. By definition,
\[
T(X_1,\dots,X_r)
=
U(X_1,\dots,X_k)\,V(X_{k+1},\dots,X_r).
\]
Its mirror bracketing is
\[
T^\ast=(V^\ast,U^\ast),
\]
so that
\[
T^\ast(Y_1,\dots,Y_r)
=
V^\ast(Y_1,\dots,Y_{r-k})\,U^\ast(Y_{r-k+1},\dots,Y_r).
\]

We prove the stated identity by induction on $r$.

For $r=1$, the statement is immediate:
\[
T(X_1)^{-1}=X_1^{-1}=T^\ast(X_1^{-1}).
\]

Assume now that $r\ge 2$ and that the result is known for all smaller numbers
of variables. Write $T=(U,V)$ as above. Then
\[
T(X_1,\dots,X_r)
=
U(X_1,\dots,X_k)\,V(X_{k+1},\dots,X_r).
\]
Since $L$ is Moufang, inversion is an anti-automorphism, so for arbitrary
subsets $P,Q\subseteq L$ one has
\[
(PQ)^{-1}=Q^{-1}P^{-1}.
\]
Applying this with
\[
P=U(X_1,\dots,X_k),\qquad Q=V(X_{k+1},\dots,X_r),
\]
we obtain
\[
T(X_1,\dots,X_r)^{-1}
=
V(X_{k+1},\dots,X_r)^{-1}\,
U(X_1,\dots,X_k)^{-1}.
\]
Now apply the induction hypothesis separately to $V$ and $U$:
\[
V(X_{k+1},\dots,X_r)^{-1}
=
V^\ast(X_r^{-1},\dots,X_{k+1}^{-1}),
\]
and
\[
U(X_1,\dots,X_k)^{-1}
=
U^\ast(X_k^{-1},\dots,X_1^{-1}).
\]
Therefore
\[
T(X_1,\dots,X_r)^{-1}
=
V^\ast(X_r^{-1},\dots,X_{k+1}^{-1})\,
U^\ast(X_k^{-1},\dots,X_1^{-1}).
\]
But this is exactly
\[
T^\ast(X_r^{-1},\dots,X_1^{-1}).
\]
This proves the first statement.

Now suppose that $A=A^{-1}$. Taking $X_1=\cdots=X_r=A$, we obtain
\[
T(A,\dots,A)^{-1}=T^\ast(A,\dots,A).
\]
Since inversion is a bijection of $L$, the map
\[
x\longmapsto x^{-1}
\]
is a bijection from $T(A,\dots,A)$ onto $T^\ast(A,\dots,A)$. Hence
\[
|T(A,\dots,A)|=|T^\ast(A,\dots,A)|.
\]
\end{proof}

\begin{remark}
For $r=4$, the five bracketings split into mirror pairs
\[
(((AA)A)A)\leftrightarrow A(A(AA)),
\qquad
((A(AA))A)\leftrightarrow A((AA)A),
\]
together with the self-mirror bracketing
\[
(AA)(AA).
\]
By Proposition~\ref{prop:mirror-bracketing}, mirror partners always have equal
cardinality when $A=A^{-1}$. Example~\ref{ex:fourfold-counterexample} shows
that this does not force equality across all bracketings.
\end{remark}

\subsection{Covering lemmas in loops}

For subsets $U,V$ of a loop $L$, define
\[
U/V:=\{u/v:u\in U,\ v\in V\},
\qquad
V\backslash U:=\{v\backslash u:v\in V,\ u\in U\}.
\]

\begin{proposition}[Loop covering lemma]
Let $L$ be a loop, and let $A,B\subseteq L$ be finite nonempty subsets. If
\[
|AB|\le K|B|,
\]
then there exists a subset $X\subseteq A$ with
\[
|X|\le K
\]
such that
\[
A\subseteq \bigcup_{x\in X}(xB)/B.
\]
\end{proposition}

\begin{proof}
Choose $X\subseteq A$ maximal subject to the sets
\[
xB\qquad (x\in X)
\]
being pairwise disjoint. Since each left translation is a bijection,
\[
|xB|=|B|
\qquad (x\in X).
\]
Hence
\[
|X|\,|B|
=
\left|\bigsqcup_{x\in X}xB\right|
\le |AB|
\le K|B|,
\]
so
\[
|X|\le K.
\]

Let $a\in A$. By maximality of $X$, the set $aB$ must intersect some $xB$ with $x\in X$. Thus there exist $b_1,b_2\in B$ such that
\[
ab_1=xb_2.
\]
By the definition of right division,
\[
a=(xb_2)/b_1\in (xB)/B.
\]
Therefore
\[
A\subseteq \bigcup_{x\in X}(xB)/B,
\]
as required.
\end{proof}

\begin{proposition}[Symmetric form in inverse-property loops]\label{prop:symmetric-ip-covering}
Let $L$ be an inverse-property loop, and let $B\subseteq L$ satisfy $B=B^{-1}$. Then for every $x\in L$,
\[
(xB)/B\subseteq (xB)B.
\]
Consequently, if $A,B\subseteq L$ are finite nonempty subsets with $B=B^{-1}$ and
\[
|AB|\le K|B|,
\]
then there exists $X\subseteq A$ with $|X|\le K$ such that
\[
A\subseteq \bigcup_{x\in X}(xB)B.
\]
\end{proposition}

\begin{proof}
In an inverse-property loop one has
\[
u/v=uv^{-1}
\]
for all $u,v\in L$, since $(uv^{-1})v=u$. Hence if $B=B^{-1}$, then
\[
(xB)/B\subseteq (xB)B.
\]
The covering statement follows immediately from the previous proposition.
\end{proof}

\begin{corollary}[Direct covering of $A^2$]
Let $L$ be a Moufang loop, and let $A\subseteq L$ be a finite $K$-approximate subloop. Then there exists a subset
\[
X\subseteq A^2
\]
with
\[
|X|\le K
\]
such that
\[
A^2\subseteq \bigcup_{x\in X}(xA)A.
\]
In particular,
\[
A^2
\]
is covered by at most $K$ subsets, each contained in $A^{\langle 4\rangle}$.
\end{corollary}

\begin{proof}
By definition,
\[
(A^2)A=(AA)A\subseteq A^{\langle 3\rangle},
\]
and therefore
\[
|(A^2)A|\le |A^{\langle 3\rangle}|\le K|A|.
\]
Apply Proposition~\ref{prop:symmetric-ip-covering} with
\[
A_{\mathrm{there}}:=A^2,\qquad B:=A.
\]
Since a Moufang loop is an inverse-property loop and $A=A^{-1}$ by definition of approximate subloop, we obtain a set
\[
X\subseteq A^2,\qquad |X|\le K,
\]
such that
\[
A^2\subseteq \bigcup_{x\in X}(xA)A.
\]
Finally, for each $x\in A^2$, every element of $(xA)A$ is represented by a fourfold product of elements of $A$ with a fixed bracketing, so
\[
(xA)A\subseteq A^{\langle 4\rangle}.
\]
\end{proof}

\begin{corollary}[Shift invariance of doubling]
Let $L$ be a commutative Moufang loop, let $x\in L$, and let $A\subseteq L$ be a finite nonempty subset. Then
\[
|(xA)^2|=|A^2|.
\]
\end{corollary}

\begin{proof}
Apply Proposition~\ref{prop:moufang-transport} with $X=Y=A$.
\end{proof}

\section{Finite-kernel reductions in Moufang loops}
\label{sec:finite-kernel}
In this section we study the basic reduction principle in arbitrary Moufang loops.

\begin{proposition}[Finite-kernel reduction]\label{prop:finite-kernel-reduction}
Let $\pi:L\twoheadrightarrow G$ be a surjective homomorphism from a loop $L$ onto a group $G$, and suppose that the kernel $N:=\Ker\pi$ is finite.
\begin{enumerate}[label=(\roman*)]
  \item If $A\subseteq L$ is a $K$-approximate subloop, then $\pi(A)$ is a $(K|N|)$-approximate subgroup of the group $G$ in the sense that
  \[
  1\in \pi(A),\qquad \pi(A)=\pi(A)^{-1},\qquad |\pi(A)^3|\le K|N|\,|\pi(A)|.
  \]
  \item If $B\subseteq G$ is an $L_0$-approximate subgroup of $G$, then $\pi^{-1}(B)$ is an $L_0$-approximate subloop of $L$.
\end{enumerate}
\end{proposition}

\begin{proof}
For~(i), the first two assertions are immediate from the homomorphism property. Since $G$ is associative, every bracketing of a triple product has the same image in $G$, and therefore
\[
\pi\bigl(A^{\langle 3\rangle}\bigr)=\pi(A)^3.
\]
Hence
\[
|\pi(A)^3|\le |A^{\langle 3\rangle}|\le K|A|.
\]
Each fiber of $\pi$ is a coset of $N$, so has cardinality $|N|$. Therefore
\[
|A|\le |N|\,|\pi(A)|,
\]
and thus
\[
|\pi(A)^3|\le K|N|\,|\pi(A)|.
\]
This proves~(i).

For~(ii), clearly $1\in \pi^{-1}(B)$ and $\pi^{-1}(B)=\pi^{-1}(B)^{-1}$. Again, since every triple product in the quotient is independent of bracketing,
\[
\pi\bigl(\pi^{-1}(B)^{\langle 3\rangle}\bigr)\subseteq B^3.
\]
Hence
\[
\pi^{-1}(B)^{\langle 3\rangle}\subseteq \pi^{-1}(B^3).
\]
Every fiber of $\pi$ has size $|N|$, so
\[
|\pi^{-1}(B)^{\langle 3\rangle}|\le |\pi^{-1}(B^3)|=|N|\,|B^3|\le L_0|N|\,|B|.
\]
But $|\pi^{-1}(B)|=|N|\,|B|$, so
\[
|\pi^{-1}(B)^{\langle 3\rangle}|\le L_0|\pi^{-1}(B)|.
\]
Thus $\pi^{-1}(B)$ is an $L_0$-approximate subloop.
\end{proof}

\begin{remark}
The proposition is purely loop-theoretic and does not use commutativity. It shows that any time a Moufang loop admits a quotient onto a genuine group with finite kernel, inverse theorems for approximate subgroups in the quotient lift immediately to the loop.
\end{remark}

\begin{proposition}[Higher product-set comparison through a finite kernel]\label{prop:higher-product-finite-kernel}
Let $\pi:L\twoheadrightarrow G$ be a surjective homomorphism from a loop $L$ onto a
group $G$, and suppose that the kernel $N:=\Ker\pi$ is finite. Let $A\subseteq L$ be a
finite subset. Then for every integer $n\ge 2$,
\[
\pi\bigl(A^{\langle n\rangle}\bigr)=\pi(A)^n.
\]
Consequently,
\[
|\pi(A)^n|\le |A^{\langle n\rangle}|\le |N|\,|\pi(A)^n|.
\]
\end{proposition}

\begin{proof}
Because $G$ is associative, every bracketing of an $n$-fold product has the same image
in $G$. Thus for every bracketing of an $n$-fold product of elements of $A$, the image of
the corresponding subset product is exactly $\pi(A)^n$. Taking the union over all
bracketings gives
\[
\pi\bigl(A^{\langle n\rangle}\bigr)=\pi(A)^n.
\]
The first inequality follows immediately. For the second, every fiber of $\pi$ has size
$|N|$, so
\[
|A^{\langle n\rangle}|\le |N|\,\bigl|\pi\bigl(A^{\langle n\rangle}\bigr)\bigr|
=|N|\,|\pi(A)^n|.
\]
\end{proof}

\begin{theorem}[Finite-kernel Ruzsa triangle inequality]
Let $\pi:L\twoheadrightarrow G$ be a surjective homomorphism from a Moufang loop $L$
onto a group $G$, and suppose that the kernel $N:=\Ker\pi$ is finite. Then for all finite
nonempty subsets $A,B,C\subseteq L$,
\[
|A|\,|BC|\le |N|\,|BA|\,|A^{-1}C|.
\]
If $A=A^{-1}$, this simplifies to
\[
|A|\,|BC|\le |N|\,|BA|\,|AC|.
\]
\end{theorem}

\begin{proof}
For each $x\in BC$, choose one pair $(b_x,c_x)\in B\times C$ such that
\[
x=b_xc_x.
\]
Define
\[
f:A\times BC\to BA\times A^{-1}C
\]
by
\[
f(a,x):=(b_xa,a^{-1}c_x).
\]
Fix $(u,v)\in BA\times A^{-1}C$. If
\[
f(a,x)=(u,v),
\]
then in the group $G$ we have
\[
\pi(u)\pi(v)
=
\pi(b_xa)\pi(a^{-1}c_x)
=
\pi(b_x)\pi(a)\pi(a)^{-1}\pi(c_x)
=
\pi(b_x)\pi(c_x)
=
\pi(x).
\]
Thus $x$ must lie in the fiber
\[
\pi^{-1}\bigl(\pi(u)\pi(v)\bigr),
\]
which has cardinality $|N|$.

Once $x$ is fixed, our choice of $(b_x,c_x)$ is fixed. Then $a$ is uniquely determined
from $u=b_xa$, namely
\[
a=b_x^{-1}u,
\]
because Moufang loops have the inverse property. Therefore every fiber of $f$ has size
at most $|N|$. Hence
\[
|A|\,|BC|
=
|A\times BC|
\le
|N|\,|BA|\,|A^{-1}C|,
\]
as claimed.
\end{proof}

\begin{corollary}[Two-generated Moufang case]
Let $M$ be a Moufang loop and let $A\subseteq M$ be a finite $K$-approximate subloop. If the generated subloop $H=\gen{A}$ can be generated by at most two elements, then $H$ is a group. Consequently $A$ is an approximate subgroup in the ordinary associative sense.
\end{corollary}

\begin{proof}
By diassociativity, every two elements of a Moufang loop generate a group~\cite{Bruck58,Pflugfelder90}. Hence if $H$ is two-generated, $H$ is a group.
\end{proof}

\begin{corollary}[Two-generated commutative Moufang case]
Let $M$ be a commutative Moufang loop and let $A\subseteq M$ be a finite
$K$-approximate subloop. If the generated subloop $H=\gen{A}$ can be generated by at
most two elements, then $H$ is an abelian group. Consequently $A$ is an approximate
subgroup in the ordinary abelian sense.
\end{corollary}

\begin{proof}
By diassociativity, every two-generated subloop of a Moufang loop is a group. Since
$H$ is also commutative, it is an abelian group.
\end{proof}

\section{Pullbacks of coset progressions in commutative Moufang loops}
\label{sec:pullbacks}
We now specialize to commutative Moufang loops. The key additional fact is that the associator subloop of a finitely generated commutative Moufang loop is finite, while the quotient by it is abelian.

Let $L$ be a commutative Moufang loop, and let
\[
\pi_L:L\to L/L'
\]
be the quotient map, where $L'=(L,L,L)$ is the associator subloop. Since
$L/L'$ is an abelian group, a natural loop-theoretic analogue in $L$ of a coset progression
is the pullback of a coset progression from the quotient.

\begin{definition}
Let $L$ be a commutative Moufang loop. A subset $P\subseteq L$ is called a
\emph{pullback of a coset progression from the quotient} if there exists a coset
progression $Q\subseteq L/L'$ such that
\[
P=\pi_L^{-1}(Q).
\]
Its rank is, by definition, the rank of $Q$.
\end{definition}

\begin{remark}
This definition is \emph{local}: it depends on the chosen commutative Moufang
loop $L$. If $H\le M$ is a subloop of a commutative Moufang loop
$M$, then such a pullback in $H$ is defined using the quotient $H/H'$, not the
ambient quotient $M/M'$. One has only
\[
H'=(H,H,H)\le M'\cap H,
\]
and in general pullbacks from $H/H'$ will not coincide with pullbacks from the image of $H$ in $M/M'$.
\end{remark}

\begin{proposition}[Canonical representatives modulo the associator subloop]
Let $H$ be a commutative Moufang loop, let $H'=(H,H,H)$, and let $\pi_H:H\to H/H'$ be the quotient map. Fix elements $x_1,\dots,x_r\in H$ and integers $\ell_1,\dots,\ell_r$. Then every word in the symbols $x_1^{\pm1},\dots,x_r^{\pm1}$ whose total exponent vector is $(\ell_1,\dots,\ell_r)$, with arbitrary ordering and arbitrary bracketing, can be written in the form
\[
\Bigl(\cdots\bigl((x_1^{\ell_1}x_2^{\ell_2})x_3^{\ell_3}\bigr)\cdots x_r^{\ell_r}\Bigr)h
\qquad\text{for some }h\in H'.
\]
\end{proposition}

\begin{proof}
Because $H$ is Moufang, powers $x_i^{\ell_i}$ are well defined. Let $W$ denote the given word and let
\[
X(\ell_1,\dots,\ell_r):=\Bigl(\cdots\bigl((x_1^{\ell_1}x_2^{\ell_2})x_3^{\ell_3}\bigr)\cdots x_r^{\ell_r}\Bigr).
\]
In the quotient group $H/H'$, both $W$ and $X(\ell_1,\dots,\ell_r)$ have the same image, namely
\[
\pi_H(x_1)^{\ell_1}\pi_H(x_2)^{\ell_2}\cdots \pi_H(x_r)^{\ell_r},
\]
because the quotient is associative and commutative. Therefore
\[
\pi_H\bigl(WX(\ell_1,\dots,\ell_r)^{-1}\bigr)=1,
\]
so
\[
WX(\ell_1,\dots,\ell_r)^{-1}\in H'.
\]
Hence
\[
W\in H'X(\ell_1,\dots,\ell_r).
\]
Since $H'$ is normal in $H$, left and right cosets coincide, and therefore
\[
H'X(\ell_1,\dots,\ell_r)=X(\ell_1,\dots,\ell_r)H'.
\]
Thus there exists $h\in H'$ such that
\[
W=X(\ell_1,\dots,\ell_r)h,
\]
as claimed.
\end{proof}

\begin{remark}
Thus, once a coset progression $Q\subseteq H/H'$ is fixed, changing lifts, order, or bracketing changes the corresponding representative in $H$ only by multiplication by an element of the finite kernel $H'$. This is a quotient-theoretic formulation of the ``associator collection'' point of view~\cite{Smith78}.
\end{remark}

\subsection{Additive-combinatorial input from abelian groups}
We shall use two standard inverse theorems in abelian groups.

\begin{theorem}[Green--Ruzsa \cite{GR07}]
For every $L\ge 1$ there exist constants $r(L),C(L)$ such that the following holds. If $B$ is a finite symmetric subset of an abelian group, with $0\in B$ and
\[
|B+B|\le L|B|,
\]
then $B$ is contained in a coset progression $Q$ of rank at most $r(L)$ and cardinality at most $C(L)|B|$.
\end{theorem}

\begin{theorem}[Gowers--Green--Manners--Tao, bounded torsion case]
For every integer $m\ge 2$ there exists $C_m>0$ such that the following holds.
If $B$ is a finite symmetric subset of an abelian group $G$ of torsion dividing $m$, with $0\in B$ and
\[
|B+B|\le L|B|,
\]
then there exist a subgroup $U\le G$ and a set $X\subseteq G$ such that
\[
|X|\le (2L)^{C_m},
\qquad
B\subseteq X+U,
\qquad
|U|\le |B|.
\]
This is the bounded-torsion case of Marton's conjecture~\cite{GGMT25, GGMT26}.
\end{theorem}

\begin{remark}
We use only these qualitative forms. No sharper quantitative form is needed for the present arguments. For each fixed $m$, the bounded-torsion theorem gives polynomial covering with exponent depending on $m$; in the exponent-$3$ application, this specializes to an absolute constant $C=C_3$.
\end{remark}

\begin{remark}
When applying inverse theorems in the abelian quotient $H/H'$, we identify multiplicative notation with additive notation. Thus statements written as $|B+B|\le L|B|$ are applied to multiplicative sets via the correspondence $B+B \leftrightarrow B^2$.
\end{remark}

\begin{theorem}[Local correspondence in a commutative Moufang loop]\label{thm:local-correspondence}
Let $M$ be a commutative Moufang loop, let $A\subseteq M$ be finite, and put
\[
H:=\gen{A},
\qquad
H':=(H,H,H),
\qquad
\pi_H:H\to H/H'.
\]
Then $H'$ is finite and $H/H'$ is an abelian group.

If $A$ is a $K$-approximate subloop of $M$, then $\pi_H(A)$ is a $(K|H'|)$-approximate subgroup of $H/H'$.

Conversely, if $B\subseteq H/H'$ is an $L$-approximate subgroup, then $\pi_H^{-1}(B)$ is an $L$-approximate subloop of $H$.
\end{theorem}

\begin{proof}
Because $A$ is finite, the subloop $H=\gen{A}$ is finitely generated. By the classical structure theory of commutative Moufang loops, $H'$ is finite and $H/H'$ is an abelian group. The two approximation statements are exactly Proposition~\ref{prop:finite-kernel-reduction} applied to the quotient map $\pi_H$.
\end{proof}

\begin{corollary}[Local covering of $A^2$ modulo the associator subloop]
With notation as above, if $A$ is a finite $K$-approximate subloop, then there exists a set
\[
X\subseteq A^2
\]
with
\[
|X|\le K|H'|
\]
such that
\[
A^2\subseteq \bigcup_{x\in X}(xA^2)H'.
\]
Equivalently, every element $u\in A^2$ can be written in the form
\[
u=(xv)h
\qquad\text{with }x\in X,\ v\in A^2,\ h\in H'.
\]
\end{corollary}

\begin{proof}
Let $B:=\pi_H(A)\subseteq H/H'$. By the local correspondence theorem, $B$ is a $(K|H'|)$-approximate subgroup of the abelian group $H/H'$. Since
\[
|B\cdot B^2|=|B^3|\le K|H'|\,|B|,
\]
the Ruzsa covering lemma in $H/H'$~\cite{Ruzsa99} yields a set
\[
\overline X\subseteq B^2
\]
with
\[
|\overline X|\le \frac{|B^3|}{|B|}\le K|H'|
\]
such that
\[
B^2\subseteq \overline X\,BB^{-1}.
\]
Because $B=B^{-1}$, this becomes
\[
B^2\subseteq \overline X\,B^2.
\]

Since $\overline X\subseteq B^2=\pi_H(A)^2=\pi_H(A^2)$, we may choose a lift $X\subseteq A^2$ of $\overline X$. If $u\in A^2$, then
\[
\pi_H(u)\in B^2\subseteq \overline X\,B^2,
\]
so there exist $x\in X$ and $v\in A^2$ such that
\[
\pi_H(u)=\pi_H(x)\pi_H(v)=\pi_H(xv).
\]
Hence
\[
u\in (xv)H'\subseteq (xA^2)H'.
\]
Since $u\in A^2$ was arbitrary, the result follows.
\end{proof}

\begin{corollary}[Local Ruzsa triangle inequality with finite-kernel loss]
Let $M$ be a commutative Moufang loop, and let $A,B,C\subseteq M$ be finite nonempty subsets. Put
\[
H:=\gen{A\cup B\cup C},
\qquad
H':=(H,H,H),
\qquad
\pi_H:H\to H/H'.
\]
Then
\[
|A|\,|BC|\le |H'|\,|AB|\,|A^{-1}C|.
\]
If $A=A^{-1}$, this simplifies to
\[
|A|\,|BC|\le |H'|\,|AB|\,|AC|.
\]
\end{corollary}

\begin{proof}
Apply the finite-kernel Ruzsa triangle inequality to the quotient map
\[
\pi_H:H\twoheadrightarrow H/H',
\]
whose kernel is $H'$. This gives
\[
|A|\,|BC|\le |H'|\,|BA|\,|A^{-1}C|.
\]
Since $M$ is commutative, one has $BA=AB$. This proves the first inequality, and the
second follows immediately when $A=A^{-1}$.
\end{proof}

\begin{corollary}[Local higher product-set comparison]
Let $M$ be a commutative Moufang loop, let $A\subseteq M$ be finite, and put
\[
H:=\gen{A},\qquad H':=(H,H,H),\qquad \pi_H:H\to H/H'.
\]
Then for every $n\ge 2$,
\[
\pi_H\bigl(A^{\langle n\rangle}\bigr)=\pi_H(A)^n
\]
and
\[
|\pi_H(A)^n|\le |A^{\langle n\rangle}|\le |H'|\,|\pi_H(A)^n|.
\]
\end{corollary}

\begin{proof}
Apply Proposition~\ref{prop:higher-product-finite-kernel} to the quotient map $\pi_H$.
\end{proof}

\section{A local structure theorem for commutative Moufang loops}
\label{sec:local-structure}
The ambient finite-generation hypothesis is unnecessary: every finite set already lies in a finitely generated subloop.

\begin{theorem}[Local structure theorem for commutative Moufang loops]\label{thm:local-structure}
Let $M$ be a commutative Moufang loop, and let $A\subseteq M$ be a finite $K$-approximate subloop. Set
\[
H:=\gen{A},\qquad H':=(H,H,H),\qquad \pi_H:H\to H/H'.
\]
Then:
\begin{enumerate}[label=(\roman*)]
  \item $H$ is finitely generated;
  \item $H'$ is finite;
  \item the image $B:=\pi_H(A)$ satisfies
  \[
  1\in B,\qquad B=B^{-1},\qquad |B^2|\le K|H'|\,|B|;
  \]
  \item $A$ is contained in the pullback $\pi_H^{-1}(Q)$ of a coset progression $Q\subseteq H/H'$, whose rank and relative size are bounded solely in terms of $K|H'|$.
\end{enumerate}
\end{theorem}

\begin{proof}
Since $A$ is finite, the subloop $H=\gen{A}$ is finitely generated. This proves~(i). By the classical structure theory of commutative Moufang loops, the associator subloop $H'$ is finite. This proves~(ii).

Now $H/H'$ is an abelian group. By Theorem~\ref{thm:local-correspondence}, the set $B=\pi_H(A)$ satisfies
\[
|B^3|\le K|H'|\,|B|.
\]
Since $1\in B$, one has $B^2\subseteq B^3$, and therefore
\[
|B^2|\le |B^3|\le K|H'|\,|B|.
\]
This proves~(iii).

Apply the Green--Ruzsa theorem in the abelian group $H/H'$ to $B$. There exists a coset progression $Q\subseteq H/H'$ containing $B$, whose rank and cardinality are bounded in terms of
\[
L:=K|H'|.
\]
Pulling back under $\pi_H$, we obtain
\[
A\subseteq \pi_H^{-1}(B)\subseteq \pi_H^{-1}(Q).
\]
Thus $A$ is contained in the pullback of a coset progression from the local quotient $H/H'$.

Finally, if $|Q|\le C(L)|B|$, then
\[
|\pi_H^{-1}(Q)|=|H'|\,|Q|\le |H'|\,C(L)|B|\le |H'|\,C(L)|A|.
\]
Since $K\ge 1$, one has $|H'|\le L$, so the relative size of the pullback is bounded solely in terms of $L=K|H'|$.
\end{proof}

\begin{remark}
The theorem is local in two ways: the ambient loop $M$ need not be finitely
generated, and the set controlling $A$ is defined using the local quotient $H/H'$
of the finitely generated subloop $H=\gen{A}$.
\end{remark}

\begin{remark}
Via the canonical representative proposition, one may view the pullback
$\pi_H^{-1}(Q)$ as a family of explicitly bracketed progression words together
with a bounded error term in the finite associator subloop $H'$.
\end{remark}

\begin{definition}
Let $F_m^{\cml}$ denote the free commutative Moufang loop on $m$ generators, and define
\[
\betafun(m):=\bigl|(F_m^{\cml})'\bigr|.
\]
By the finiteness statement for associator subloops in finitely generated commutative
Moufang loops, this is finite for every $m\in\mathbb{N}$~\cite{Bruck58}.
The quantity $\betafun(m)$ is related to classical questions on the structure and
enumeration of free commutative Moufang loops; see, for example, \cite{Smith82}.
\end{definition}

\begin{theorem}[Uniform \texorpdfstring{$m$}{m}-generated theorem]
Fix $m\in\mathbb{N}$. Let $M$ be a commutative Moufang loop, and let $A\subseteq M$ be a finite $K$-approximate subloop. Suppose that $H=\gen{A}$ can be generated by at most $m$ elements. Then $A$ is contained in the pullback of a coset progression from the local quotient $H/H'$ whose rank and relative size are bounded solely in terms of $K$ and $m$.
\end{theorem}

\begin{proof}
Let $H=\gen{A}$. Since $H$ can be generated by at most $m$ elements, there exists a surjective homomorphism
\[
\varphi:F_m^{\cml}\twoheadrightarrow H.
\]
Associators are preserved by homomorphisms, so $\varphi\bigl((F_m^{\cml})'\bigr)=H'$. Therefore $H'$ is a homomorphic image of $(F_m^{\cml})'$, and hence
\[
|H'|\le \betafun(m).
\]
Now apply Theorem~\ref{thm:local-structure}. The pullback of a coset progression containing $A$ has parameters depending only on $K|H'|$, and therefore only on $K\betafun(m)$, hence only on $K$ and $m$.
\end{proof}

\begin{theorem}[No elements of order \texorpdfstring{$3$}{3}]\label{thm:no-order-three}
Let $M$ be a commutative Moufang loop. Suppose that $M$ contains no nonidentity
element of order $3$. Then $M$ is an abelian group.
\end{theorem}

\begin{proof}
Let $M'=(M,M,M)$ be the associator subloop. By the classical structure theory of
commutative Moufang loops, the quotient $M/M'$ is an abelian group. On the other
hand, every element of $M'$ has cube equal to $1$. Indeed, in a commutative Moufang
loop the cube map is a centralizing endomorphism and is trivial on the associator
subloop. Therefore every element of $M'$ has order dividing $3$.

By hypothesis, $M$ has no nonidentity element of order $3$, so $M'=\{1\}$. Hence
$M\cong M/M'$ is an abelian group.
\end{proof}

\begin{corollary}
Let $M$ be a commutative Moufang loop, let $A\subseteq M$ be a finite
$K$-approximate subloop, and put $H=\gen{A}$. If $H$ contains no nonidentity element
of order $3$, then $H$ is an abelian group. Consequently $A$ is an approximate
subgroup in the ordinary abelian sense.
\end{corollary}

\begin{proof}
Apply Theorem~\ref{thm:no-order-three} to the commutative Moufang loop $H$.
\end{proof}

\section{Bounded torsion in the local quotient}
\label{sec:bounded-torsion}
We note the stronger covering consequence available when the local abelian quotient has bounded torsion.

\begin{definition}
A loop $L$ has \emph{exponent dividing $m$} if $x^m=1$ for all $x\in L$.
\end{definition}

\begin{theorem}[Local polynomial covering in bounded torsion]\label{thm:local-bounded-torsion-covering}
Let $M$ be a commutative Moufang loop, let $A\subseteq M$ be a finite
$K$-approximate subloop, and put
\[
H:=\gen{A},\qquad H':=(H,H,H),\qquad \pi_H:H\to H/H'.
\]
Assume that the abelian group $H/H'$ has torsion dividing $m$.
Then there exists a constant $C_m>0$ and a finite subloop $J\le H$ such that:
\begin{enumerate}[label=(\roman*)]
  \item $A$ is covered by at most $(2K|H'|)^{C_m}$ left cosets of $J$;
  \item $|J|\le |H'|\,|A|$.
\end{enumerate}
In particular, if $M$ has exponent dividing $m$, then the conclusion holds.
\end{theorem}

\begin{proof}
Let $B:=\pi_H(A)\subseteq H/H'$. By Theorem~\ref{thm:local-correspondence},
\[
|B^3|\le K|H'|\,|B|.
\]
Since $1\in B$, it follows that $B^2\subseteq B^3$, so
\[
|B^2|\le K|H'|\,|B|.
\]
Set
\[
L:=K|H'|.
\]
Applying the bounded-torsion Marton theorem to $B$ in the abelian group $H/H'$, we obtain a subgroup $U\le H/H'$ and a set $\overline X$ such that
\[
|\overline X|\le (2L)^{C_m},
\qquad
B\subseteq \overline X+U,
\qquad
|U|\le |B|,
\]
where $C_m$ depends only on $m$.

Define
\[
J:=\pi_H^{-1}(U)\le H.
\]
Since the kernel of $\pi_H$ is $H'$, each fiber has size $|H'|$, and therefore
\[
|J|=|H'|\,|U|\le |H'|\,|B|\le |H'|\,|A|.
\]
Choose a set $X\subseteq H$ of lifts of the elements of $\overline X$. Then
\[
A\subseteq \pi_H^{-1}(B)\subseteq XJ.
\]
Thus $A$ is covered by at most
\[
|X|=|\overline X|\le (2L)^{C_m}=(2K|H'|)^{C_m}
\]
left cosets of $J$.

Finally, if $M$ has exponent dividing $m$, then so does the quotient $H/H'$, which therefore has torsion dividing $m$. This proves the final assertion.
\end{proof}

\begin{corollary}[Polynomial covering by a finite subloop in exponent \texorpdfstring{$3$}{3}]
Let $M$ be a commutative Moufang loop of exponent $3$, and let
$A\subseteq M$ be a finite $K$-approximate subloop. Put $H=\gen{A}$ and let
$d=d(H)$ be the minimal size of a generating set of $H$. Then there exists an
absolute constant $C>0$ and a finite subloop $J\le H$ such that:
\begin{enumerate}[label=(\roman*)]
  \item $A$ is covered by at most $(2K\betafun(d))^C$ left cosets of $J$;
  \item $|J|\le \betafun(d)|A|$.
\end{enumerate}
\end{corollary}

\begin{proof}
By definition of $\betafun(d)$, one has $|H'|\le \betafun(d)$. Apply Theorem~\ref{thm:local-bounded-torsion-covering} with $m=3$.
\end{proof}

\begin{remark}
The bounded-torsion conclusion is stronger than mere containment in a pullback of a coset progression:
one obtains control by a finite subloop, with polynomial covering number in
the effective local parameter $K|H'|$. The dependence is not uniform in the torsion
parameter $m$, since the exponent $C_m$ comes from the bounded-torsion Marton theorem.
\end{remark}

\section*{Acknowledgements}
The author wishes to thank Emmanuel Breuillard for a number of helpful
discussions and advice on this topic during the time of his PhD thesis.

\end{document}